\documentclass[leqno,letterpaper]{etna}
\usepackage[sorting=none]{biblatex} %Imports biblatex package

\usepackage{amsmath}
\usepackage{graphicx}
\usepackage{amssymb}
\usepackage{algorithm}
\usepackage{algorithmic}
\usepackage{multirow}

\title{Augmented unprojected Krylov subspace methods%\thanks{%
}

\author{Liam Burke\footnotemark[2]
        \and Kirk M. Soodhalter\footnotemark[3]}

\shorttitle{Augmented unprojected Krylov subspace methods} 
\shortauthor{L. Burke and K. M. Soodhalter}

\begin{document}

\maketitle

\renewcommand{\thefootnote}{\fnsymbol{footnote}}

\footnotetext[1]{School Of Mathematics, Trinity College Dublin, College Green, Dublin 2,  \url{burkel8@tcd.ie}}
\footnotetext[2]{School Of Mathematics, Trinity College Dublin, College Green, Dublin 2,  \url{ksoodha@maths.tcd.ie}}

\begin{abstract}
Augmented Krylov subspace methods aid in accelerating the convergence of a standard Krylov subspace method by including additional vectors in the search space. A residual projection framework based on residual (Petrov-) Galerkin constraints was presented in [Gaul et al. SIAM J. Matrix Anal.
Appl 2013], and later generalised in a recent survey on subspace recycling iterative methods [Soodhalter et al. GAMM-Mitt. 2020]. The framework describes augmented Krylov subspace methods in terms of applying a standard Krylov subspace method to an appropriately projected problem.

In this work we show that the projected problem has an equivalent unprojected formulation, and that viewing the framework in this way provides a similar description for the class of \textit{unprojected} augmented Krylov subspace methods. We derive the first unprojected augmented Full Orthogonalization Method (FOM), and demonstrate its effectiveness as a recycling method. We then show how the R$^{3}$GMRES algorithm fits within the framework. We show that unprojected augmented short recurrence methods fit within the framework, but can only be implemented in practice under certain conditions on the augmentation subspace. We demonstrate this using the Augmented Conjugate Gradient (AugCG) algorithm as an example.
\end{abstract}



\section{Introduction}
It is well known that iterative Krylov subspace methods are the most practical and efficient means of approximating the solution to a linear system of equations with large and sparse coefficient matrix. In High-Performance Computing (HPC) applications, it is of interest to accelerate the convergence of these methods and reduce the number of matrix vector operations required while still maintaining approximation accuracy within the specified level of tolerance. One approach is to use \textit{augmented} Krylov subspace methods
wherein the solution is approximated over the sum of an iteratively built $j$-dimensional Krylov subspace $\mathcal{V}_{j}$ and some other fixed $k$-dimensional subspace $\mathcal{U}$ known as the \textit{augmentation subspace}.  The subspace $\mathcal{U}$ is carefully chosen  to contain useful information to aid in accelerating the convergence of the iterative solver, and is commonly taken to be the subspace spanned by approximate eigenvectors corresponding to the $k$ smallest eigenvalues of the coefficient matrix. The technique of \textit{Krylov subspace recycling} is a special case of an augmented Krylov subspace methods wherein the subspace $\mathcal{U}$ comes from the solution of a previous system solve. There is a large body of literature dedicated to the study of augmentation and recycling (see eg. \cite{gaul2013framework, parks2006recycling,gutknecht2014deflated, gaul2014recycling}). In particular, we refer the reader to the recent survey \cite{Soodhalter2020ASO} and references therein which provide a broad overview of the topic.

Augmented Krylov subspace methods are typically described in terms of a well known residual projection framework which describes these methods in terms of applying a standard Krylov subspace method to some projected problem, followed by taking an additional solution correction from the augmentation subspace. The framework was originally proposed in terms of typical residual constraints in works such as \cite{gaul2014recycling, gaul2013framework, gutknecht2012spectral,  gutknecht2014deflated}, but the survey \cite{Soodhalter2020ASO} extends the framework to the case of arbitrary search and constraint spaces.

Despite the original view of this framework, not all augmented Krylov subspace methods work using a Krylov subspace built from a projected matrix. It is often of interest to use Krylov subspaces built from the original unprojected matrix, and we refer to such methods as \textit{unprojected} augmented Krylov subspace methods. Unprojected methods are often of interest in the solution of ill-posed inverse problems (see e.g, \cite{dong2014r3gmres}), or in situations where the projector is too expensive to apply or there is a strict constraint on the amount of computations that can be performed per iteration (see e.g, \cite{erhel2000augmented, ramlau2021augmented}).

The main goal of this paper is to show that the projected problem arising in the framework  \cite{gaul2013framework,Soodhalter2020ASO} has an equivalent unprojected formulation, and that by viewing this projected problem in its unprojected form, the framework can be extended to describe the class of unprojected augmented Krylov subspace methods. We observe the difference in describing a projected method and an unprojected method within this framework is entirely dependent on whether the projected or unprojected problem is solved and that the the final correction from the augmentation subspace  always remains the same.

We demonstrate how to easily derive an unprojected augmented Krylov subspace method from this view of the framework simply by making appropriate choices of subspaces. For non-Hermitian matrices, we develop the first unprojected augmented Full Orthogonalization Method (FOM) and demonstrate its effectiveness as a recycling method. We then
show that for choices of constraint subspaces corresponding to residual minimizing methods, the R$^{3}$GMRES \cite{dong2014r3gmres} algorithm arises naturally from this framework.

Finally, for a Symmetric Positive Definite (SPD) matrix, in the projected formulation, the SPD properties of the matrix are typically inherited into the projected problem as shown in \cite{Soodhalter2020ASO}. We use the framework to show that an analogy of this result generally does not hold for unprojected methods unless additional constraints are imposed on the augmentation subspace. While this leads to practical limitations for the development of unprojected augmented short recurrence methods, it is still possible to develop such methods from the framework, and indeed the framework allows us to understand which choices of augmentation subspace yield practical methods. We demonstrate this by using the framework to derive the Augmented Conjugate Gradient algorithm (AugCG) \cite{erhel2000augmented}.
\section{An overview of Krylov subspace methods}\label{sec:krylov_overview}
We shortly summarize the basics of a residual projection framework for solving linear systems of the form
\begin{equation} \label{linearsystem}
    \textbf{A}\textbf{x} = \textbf{b}
\end{equation}
where $\textbf{A} \in \mathbb{C}^{n \times n}$ , $\textbf{b} \in \mathbb{C}^{n}$, and the matrix $\textbf{A}$ is assumed to be large and sparse.
A projection method projects the problem  (\ref{linearsystem}) onto the iteratively built $j$ dimensional subspace $\mathcal{V}_{j}$ and solves a smaller representation of the problem on $\mathcal{V}_{j}$ using a direct method (see e.g, \cite{saad2003iterative}). Assuming knowledge of an initial approximation $\textbf{x}_{0}$ and associated residual $\textbf{r}_{0}$, we seek a solution correction $\textbf{t}_{j} \in \mathcal{V}_{j}$ such that the updated residual $\textbf{r}_{j}$ satisfies the orthogonality constraint
\begin{equation}\label{eq:residual_constraint}
\textbf{r}_{j} = \textbf{b} - \textbf{A}(\textbf{x}_{0} + \textbf{t}_{j}) \perp \widetilde{\mathcal{V}}_{j},
\end{equation}
where $\widetilde{\mathcal{V}}_{j}$ is some other $j$ dimensional subspace.

Since $\textbf{t}_{j} \in \mathcal{V}_{j}$ it can be written as $\textbf{t}_{j} = \textbf{V}_{j} \textbf{y}_{j}$ where $\textbf{V}_{j}$ is a matrix with columns forming a basis for $\mathcal{V}_{j}$ and $\textbf{y}_{j} \in \mathbb{C}^{j}$ is a coefficients vector. Application of the residual orthogonality constraint (\ref{eq:residual_constraint}) yields an equation for $\textbf{y}_{j}$ given by
\begin{equation}\label{eq:probfory}
 (\widetilde{\textbf{V}}_{j}^{*} \textbf{A} \textbf{V}_{j}) \textbf{y}_{j} = \widetilde{\textbf{V}}_{j}^{*} \textbf{r}_{0},
\end{equation}
where $\widetilde{\textbf{V}}_{j} \in \mathbb{C}^{n \times j}$ is a matrix with columns forming a basis for $\widetilde{\mathcal{V}}_{j}$.
After solving (\ref{eq:probfory}) for $\textbf{y}_{j}$, the updated solution approximation $\textbf{x}_{j}$ can be computed via $\textbf{x}_{j} = \textbf{x}_{0} + \textbf{V}_{j} \textbf{y}_{j}$.

A practical implementation of a projection method requires an appropriate choice of subspaces $\mathcal{V}_{j}$ and $\widetilde{\mathcal{V}}_{j}$.
Krylov subspace methods are the class of methods which choose $\mathcal{V}_{j}$ to be the $j^{th}$ order \textit{Krylov subspace} generated from the matrix $\textbf{A}$ and vector $\textbf{r}_{0}$ defined by
\[
    \mathcal{K}_{j}(\textbf{A},\textbf{r}_{0}) = \mbox{span}\{ \textbf{r}_{0}, \textbf{A}\textbf{r}_{0}, \textbf{A}^{2} \textbf{r}_{0}, \ldots , \textbf{A}^{j-1} \textbf{r}_{0} \}.
\]

For non-Hermitian matrices the Arnoldi process can be used to construct an orthonormal basis for $\mathcal{K}_{j}(\textbf{A},\textbf{b})$ stored as columns in the matrix $\textbf{V}_{j+1} \in \mathbb{C}^{n \times (j+1)}$, and an upper Hessenberg matrix $\overline{\textbf{H}}_{j} \in \mathbb{C}^{(j+1) \times j}$ which stores the orthogonalization coefficients.
At the end of the $j$-length cycle of the Arnoldi process the following \textit{Arnoldi relation} is satisfied
\begin{equation}\label{eq:arnoldi}
    \textbf{A} \textbf{V}_{j} = \textbf{V}_{j+1} \overline{\textbf{H}}_{j} = \textbf{V}_{j} \textbf{H}_{j} + h_{j+1,j}\textbf{v}_{j+1} \textbf{e}^{T}_{j},
\end{equation}
with $\textbf{H}_{j} \in \mathbb{C}^{j \times j}$ the matrix $\overline{\textbf{H}}_{j}$ with last row deleted, and $h_{j+1,j}$  the entry at row $j+1$ and column $j$ of  $\overline{\textbf{H}}_{j}$. 

Different Krylov subspace methods are distinguished by the choice of constraint space $\widetilde{\mathcal{V}}_{j}$. The Full Orthogonalization Method (FOM) \cite{saad1981krylov} takes $\widetilde{\mathcal{V}}_{j} = \mathcal{V}_{j} = \mathcal{K}_{j}(\textbf{A},\textbf{b})$, leading to the following  $j \times j$ linear system for $\textbf{y}_{j}$
\begin{equation}\label{eq:FOM}
       \textbf{H}_{j} \textbf{y}_{j} = \| \textbf{r}_{0} \| \textbf{e}_{1}.
\end{equation}
The Generalized Minimum Residual Method (GMRES) method \cite{saad1986gmres} takes $\widetilde{\mathcal{V}}_{j} = \textbf{A} \mathcal{K}_{j}(\textbf{A},\textbf{b})$, leading to the following rectangular least squares problem for $\textbf{y}_{j}$
\begin{equation}\label{eq:GMRES}
    \textbf{y}_{j} = \operatorname*{argmin}_{\textbf{z} \in \mathbb{C}^{j}}  \left\| \| \textbf{r}_{0} \| \textbf{e}_{1} - \overline{\textbf{H}}_{j} \textbf{z} \right\|.
\end{equation}
The restarted FOM and GMRES algorithm are briefly outlined in Algorithm \ref{alg:restartedFOMandGMRES}.
\begin{algorithm}[H]
\caption{Restarted FOM and GMRES
\label{alg:restartedFOMandGMRES}} 
\begin{algorithmic}[1]
\STATE{\textbf{Input:} $\textbf{A} \in \mathbb{C}^{n \times n}$,  $\textbf{b} \in \mathbb{C}^{n}$, $\textbf{x} \in \mathbb{C}^{n}$}
\STATE{$\textbf{r} = \textbf{b} - \textbf{A} \textbf{x}$}
\WHILE{not converged}
\STATE{Build a basis for the Krylov subspace $\mathcal{K}_{j}(\textbf{A},\textbf{r})$ via the Arnoldi algorithm generating $\textbf{V}_{j+1} \in \mathbb{C}^{n \times (j+1)}$ and $\overline{\textbf{H}}_{j} \in \mathbb{C}^{(j+1) \times j}$.}
  \STATE{Compute coefficients $\textbf{y}_{j}$ using (\ref{eq:FOM}) for FOM or (\ref{eq:GMRES}) for GMRES}
  \STATE{$\textbf{x}  \leftarrow \textbf{x} + \textbf{V}_{j}\textbf{y}_{j}$}
  \STATE{$\textbf{r} \leftarrow \textbf{r} - \textbf{V}_{j+1} \overline{\textbf{H}}_{j} \textbf{y}_{j}$}
\ENDWHILE
\end{algorithmic}
\end{algorithm}

For Hermitian matrices $\textbf{A}$, the Lanczos algorithm can be used in place of Arnoldi. The resulting Hessenberg matrix in (\ref{eq:arnoldi}) is symmetric tridiagonal and often denoted as $\textbf{T}_{j}$. The FOM solution update is given by
\begin{align}\label{eq:fom_solution}
    \textbf{x}_{j} = \textbf{x}_{0} + \textbf{V}_{j} \textbf{T}_{j}^{-1} (\| \textbf{r}_{0} \| \textbf{e}_{1}),
\end{align}
and the associated residuals remain orthogonal. 
If $\textbf{A}$ is also positive definite, it omits an LU factorization $\textbf{T}_{j} = \textbf{L}_{j} \textbf{W}_{j}$
where $\textbf{L}_{j}$ is is unit lower bidiagonal and $\textbf{W}_{j}$ is upper bidiagonal. 
If we define $\textbf{P}_{j} = \textbf{V}_{j} \textbf{W}_{j}^{-1}$ and $\textbf{z}_{j} = \textbf{L}_{j}^{-1}(\| \textbf{r}_{0} \| \textbf{e}_{1})$ then (\ref{eq:fom_solution}) can be written as
\[
    \textbf{x}_{j} = \textbf{x}_{0} + \textbf{P}_{j} \textbf{z}_{j}.
\]
Due to the structure of $\textbf{W}_{j}$, the columns of $\textbf{P}_{j}$ can be constructed progressively by
\[
    \textbf{p}_{j} = (\textbf{v}_{j} - \beta_{j} \textbf{p}_{j-1})/\eta_{j},
\]
where $\eta_{j}$ is the $j^{th}$ entry on the main diagonal of $\textbf{W}_{j}$ and $\beta_{j}$ is a scalar produced by the Lanczos process. The structure of $\textbf{L}_{j}$ then allows $\textbf{z}_{j}$ to be written as
\[
    \textbf{z}_{j} = \begin{bmatrix}
          \textbf{z}_{j-1} \\
           \xi_{j}
    \end{bmatrix},  \hspace{1.1cm} \xi_{j} \in \mathbb{C},
\]
and the full solution approximation is then given by
\[
    \textbf{x}_{j} = \textbf{x}_{j-1} + \xi_{j} \textbf{p}_{j}.
\]
It can be shown (see e.g, \cite{erhel2000augmented, saad2003iterative}) that the vectors $\textbf{p}_{j}$ form an $\textbf{A}$-conjugate set. Explicitly imposing this conjugacy condition along with the orthogonality of the residual vectors yields the updates
\[
    \textbf{x}_{j} = \textbf{x}_{j-1} + \alpha_{j} \textbf{p}_{j},
\]
with $\alpha_{j} = \textbf{r}_{j}^{T} \textbf{r}_{j}/ \textbf{r}_{j}^{T}\textbf{A}\textbf{p}_{j}$ and
\[
\textbf{p}_{j+1} = \textbf{r}_{j+1} + \beta_{j} \textbf{p}_{j},
\]
with $\beta_{j} = -\textbf{p}_{j}^{T} \textbf{A} \textbf{r}_{j+1}/\textbf{p}_{j}^{T}\textbf{A}\textbf{p}_{j}$.

The solution updates no longer require storage of the full Arnoldi basis. This method is known as the Conjugate Gradient (CG) algorithm \cite{hestenes1952methods} and is outlined in Algorithm \ref{alg:CG}.

\begin{algorithm}[H]
\caption{Conjugate Gradient (CG) \label{alg:CG}}
\begin{algorithmic}[1]
\STATE{\textbf{Input:} $\textbf{A} \in \mathbb{C}^{n \times n}$ , $ \textbf{b} \in \mathbb{C}^{n}$, initial residual $\textbf{r}_{0} \in \mathbb{C}^{n}$ from initial approximation $\textbf{x}_{0} \in \mathbb{C}^{n}$}
\FOR{$j = 0,1, \hdots$ until convergence}
\STATE{$\alpha_{j} = \textbf{r}_{j}^{T} \textbf{r}_{j} / \textbf{p}_{j}^{T} \textbf{A} \textbf{p}_{j}$}
\STATE{$\textbf{x}_{j+1} = \textbf{x}_{j} + \alpha_{j} \textbf{p}_{j}$}
\STATE{$\textbf{r}_{j+1} = \textbf{r}_{j} - \alpha_{j} \textbf{A} \textbf{p}_{j}$}
\STATE{$\beta_{j} = \textbf{r}_{j+1}^{*} \textbf{r}_{j+1}/\textbf{r}_{j}^{*}\textbf{r}_{j}$}
\STATE{$\textbf{p}_{j+1} = \textbf{r}_{j+1} + \beta_{j} \textbf{p}_{j}$}
\ENDFOR
\end{algorithmic}
\end{algorithm}

Throughout the remainder of the paper we use the term \textit{standard Krylov subspace method} to refer to methods which do not use augmentation. The FOM, GMRES and CG methods outlined in this section are such examples. 
For any augmented Krylov subspace method, we refer to its \textit{corresponding standard Krylov subspace method} as the original non-augmented version of the method. For example, the corresponding non-augmented version of augmented GMRES is just standard GMRES.  

\section{A framework for augmented Krylov subspace methods}\label{sec:framework}
In this section we outline the augmented Krylov subspace framework from works such as \cite{ gaul2013framework, gutknecht2014deflated, gaul2014recycling,gutknecht2012spectral} for solving (\ref{linearsystem}). We describe the framework in its most most general form as presented in the survey \cite{Soodhalter2020ASO}.

Beginning with an an initial approximation $\textbf{x}_{0} \in \mathbb{C}^{n}$, an augmented Krylov subspace method for solving (\ref{linearsystem}) computes a solution correction $\textbf{s}_{j} \in \mathbb{C}^{n}$ from the space $\mathcal{U}$, and a correction $\textbf{t}_{j} \in \mathbb{C}^{n}$ from the Krylov space $\mathcal{V}_{j}$. Let  $\textbf{U} \in \mathbb{C}^{n \times k}$ and $\textbf{V}_{j} \in \mathbb{C}^{n \times j}$ be matrices with columns forming a basis for $\mathcal{U}$ and $\mathcal{V}_{j}$ respectively, then we take the updated solution approximation $\textbf{x}_{j}$ at step $j$ to be
\[
    \textbf{x}_{j} = \textbf{x}_{0} + \textbf{s}_{j} + \textbf{t}_{j} = \textbf{x}_{0} + \textbf{U} \textbf{z}_{j} +  \textbf{V}_{j}\textbf{y}_{j},
\]
where $\textbf{z}_{j} \in \mathbb{C}^{k}$ and $\textbf{y}_{j} \in \mathbb{C}^{j}$ are coefficient vectors determined by imposing the following orthogonality condition
\begin{equation}\label{eq:framework_residual_constrict}
    \textbf{r}_{j} = \textbf{b} - \textbf{A}(\textbf{x}_{0} + \textbf{U} \textbf{z}_{j} + \textbf{V}_{j} \textbf{y}_{j}) \perp \widetilde{\mathcal{U}} + \widetilde{\mathcal{V}}_{j},
\end{equation}
where $\widetilde{\mathcal{U}}$ is some 
$k$ dimensional subspace (not necessarily the same as $\mathcal{U}$), and $\widetilde{\mathcal{V}}_{j}$ is some Krylov subspace (not necessarily the same as $\mathcal{V}_{j}$). Application of the residual constraint (\ref{eq:framework_residual_constrict}) leads to a $(k + j) \times (k + j)$ linear system for $\textbf{z}_{j}$ and $\textbf{y}_{j}$ given by
\begin{equation}\label{linear_system}
    \begin{bmatrix}
        \widetilde{\textbf{U}}^{*} \textbf{A} \textbf{U} &  \widetilde{\textbf{U}}^{*} \textbf{A} \textbf{V}_{j} \\\widetilde{\textbf{V}}^{*}_{j} \textbf{A} \textbf{U} &  \widetilde{\textbf{V}}^{*}_{j} \textbf{A} \textbf{V}_{j}
    \end{bmatrix} \begin{bmatrix} \textbf{z}_{j} \\ \textbf{y}_{j} \end{bmatrix} = \begin{bmatrix} 
       \widetilde{\textbf{U}}^{*} \textbf{r}_{0} \\
       \widetilde{\textbf{V}}^{*}_{j} \textbf{r}_{0}
    \end{bmatrix}.
\end{equation}

A block LU factorization of the linear system (\ref{linear_system}) decouples an equation for $\textbf{y}_{j}$ (independent of $\textbf{z}_{j}$) given as the second block row of the linear system
\begin{align*}
     \begin{bmatrix}
        \widetilde{\textbf{U}}^{*} \textbf{A} \textbf{U} &  \widetilde{\textbf{U}}^{*} \textbf{A} \textbf{V}_{j} \\ \textbf{0} &  \widetilde{\textbf{V}}^{*}_{j} \textbf{A} \textbf{V}_{j} - \widetilde{\textbf{V}}^{*}_{j} \textbf{A} \textbf{U} (\widetilde{\textbf{U}}^{*} \textbf{A} \textbf{U})^{-1} (\widetilde{\textbf{U}}^{*} \textbf{A} \textbf{V}_{j})
    \end{bmatrix} \begin{bmatrix} \textbf{z}_{j} \\ \textbf{y}_{j} \end{bmatrix} =   \begin{bmatrix}
        \widetilde{\textbf{U}}^{*} \textbf{r}_{0} \\
        \widetilde{\textbf{V}}^{*}_{j} \textbf{r}_{0} - \widetilde{\textbf{V}}^{*}_{j} \textbf{A} \textbf{U}(\widetilde{\textbf{U}}^{*} \textbf{A} \textbf{U})^{-1}\widetilde{\textbf{U}}^{*} \textbf{r}_{0}
     \end{bmatrix},
\end{align*}
which we can write as
\begin{align}\label{eq:eqfory}
    \widetilde{\textbf{V}}^{*}_{j}( \textbf{I} -\textbf{A} \textbf{U}(\widetilde{\textbf{U}}^{*} \textbf{A} \textbf{U})^{-1}\widetilde{\textbf{U}}^{*}) \textbf{A} \textbf{V}_{j} \textbf{y}_{j} = \widetilde{\textbf{V}}^{*}_{j}( \textbf{I} - \textbf{A} \textbf{U}(\widetilde{\textbf{U}}^{*} \textbf{A} \textbf{U})^{-1} \widetilde{\textbf{U}}^{*}) \textbf{r}_{0}.
\end{align}

Once the vector $\textbf{y}_{j}$ is computed, the vector $\textbf{z}_{j}$ can the be computed as
\begin{equation}\label{zcorrection}
    \textbf{z}_{j} = (\widetilde{\textbf{U}}^{*} \textbf{A} \textbf{U})^{-1} \widetilde{\textbf{U}}^{*} \textbf{r}_{0} - (\widetilde{\textbf{U}}^{*} \textbf{A} \textbf{U})^{-1} (\widetilde{\textbf{U}}^{*} \textbf{A} \textbf{V}_{j}) \textbf{y}_{j},
\end{equation}
leading to the full solution update 
\begin{align}\label{eq:framework_full_solution_update}
    \textbf{x}_{j} = \textbf{x}_{0} +  \textbf{U} (\widetilde{\textbf{U}}^{*} \textbf{A} \textbf{U})^{-1} \widetilde{\textbf{U}}^{*} \textbf{r}_{0} + \textbf{V}_{j} \textbf{y}_{j} - \textbf{U} (\widetilde{\textbf{U}}^{*} \textbf{A} \textbf{U})^{-1} (\widetilde{\textbf{U}}^{*} \textbf{A} \textbf{V}_{j}) \textbf{y}_{j},
\end{align}
and corresponding updated residual
\begin{align}\label{eq:framework_full_residual_update}
    \textbf{r}_{j} = \textbf{r}_{0} - \textbf{A} \textbf{U}(\widetilde{\textbf{U}}^{*} \textbf{A} \textbf{U})^{-1} \widetilde{\textbf{U}}^{*} \textbf{r}_{0}
     - \textbf{A} \textbf{V}_{j} \textbf{y}_{j}
    + \textbf{A} \textbf{U} (\widetilde{\textbf{U}}^{*} \textbf{A} \textbf{U})^{-1} \widetilde{\textbf{U}}^{*} \textbf{A} \textbf{V}_{j} \textbf{y}_{j}.
\end{align}

The main feature of the framework comes from defining the projector $\Phi_{\textbf{A}\mathcal{U},\widetilde{\mathcal{U}}^{\perp}} := \textbf{A} \textbf{U}(\widetilde{\textbf{U}}^{*}\textbf{A} \textbf{U})^{-1} \widetilde{\textbf{U}}^{*}$, from which we can then see from (\ref{eq:eqfory}) that computing $\textbf{y}_{j}$ is equivalent to applying the following standard Krylov subspace method.
\\ \newline
\noindent\fbox{%
    \parbox{\textwidth}{
   Find $\textbf{t}_{j} \in \mathcal{V}_{j}$  as an approximate solution to the problem 
\begin{equation}\label{eq:framework_projectedproblem}
   (\textbf{I} - \Phi_{\textbf{A} \mathcal{U},\widetilde{\mathcal{U}}^{\perp}}) \textbf{A} \textbf{t} =  (\textbf{I} - \Phi_{\textbf{A} \mathcal{U},\widetilde{\mathcal{U}}^{\perp}}) \textbf{r}_{0}
\end{equation}
 subject to $\textbf{r}_{j} \perp \widetilde{\mathcal{V}}_{j}$
    }
}
\\ \newline

In other words, if we perform the residual projection $ \widehat{\textbf{r}}  = (\textbf{I} - \Phi_{\textbf{A} \mathcal{U},\widetilde{\mathcal{U}}^{\perp}}) \textbf{r}_{0},$ and take the subspace  $\mathcal{V}_{j}$ to be the Krylov subspace built from the projected operator and  projected residual 
\[
    \mathcal{V}_{j} = \mathcal{K}_{j}((\textbf{I} - \Phi_{\textbf{A} \mathcal{U},\widetilde{\mathcal{U}}^{\perp}}) \textbf{A}, \widehat{\textbf{r}}),
\]
then the augmented Krylov subspace method requires application of the corresponding standard Krylov subspace method to an appropriately projected problem, followed by taking an additional solution correction from the augmentation subspace.

The initial residual projection corresponds to the first two terms of (\ref{eq:framework_full_residual_update}), while the resulting effect of this projection on the solution corresponds to the first two terms in the solution update (\ref{eq:framework_full_solution_update}).

In order to build a basis for the Krylov subspace 
$\mathcal{K}_{j}((\textbf{I} - \Phi_{\textbf{A} \mathcal{U},\widetilde{\mathcal{U}}^{\perp}}) \textbf{A}, \widehat{\textbf{r}})$, it is not necessary to explicitly construct the operator $(\textbf{I} - \Phi_{\textbf{A} \mathcal{U},\widetilde{\mathcal{U}}^{\perp}}) \textbf{A}$. Instead the original matrix $\textbf{A}$ is applied to each Arnoldi vector and the action of $(\textbf{I} - \Phi_{\textbf{A} \mathcal{U},\widetilde{\mathcal{U}}^{\perp}})$ can be performed implicitly throughout the Arnoldi process. For an orthogonal projection, this is done using an additional orthogonalization of the Arnoldi vectors. The coefficients of this orthogonalization are stored in the matrix $\textbf{B}_{j} = (\widetilde{\textbf{U}}^{*} \textbf{A} \textbf{U})^{-1} \widetilde{\textbf{U}}^{*} \textbf{A} \textbf{V}_{j}$, which is also built iteratively throughout the process. For an oblique projection it can be done in a similar, numerically stable way using appropriate choices for $\textbf{U}$ and $\textbf{C}$ through a procedure which closely resembles Gram Schmidt. 

Building the basis for $\mathcal{K}_{j}((\textbf{I} - \Phi_{\textbf{A} \mathcal{U},\widetilde{\mathcal{U}}^{\perp}}) \textbf{A}, \widehat{\textbf{r}})$ via the  Arnoldi process yields the following Arnoldi relation
\begin{align}\label{eq:framework_projected_arnoldi_relation}
      (\textbf{I} - \Phi_{\textbf{A} \mathcal{U},\widetilde{\mathcal{U}}^{\perp}}) \textbf{A} \textbf{V}_{j} = \textbf{V}_{j+1} \overline{\textbf{H}}_{j}
\end{align}      
which we can rewrite as 
\[
 \textbf{A}\textbf{V}_{j} = \textbf{A} \textbf{U} \textbf{B}_{j} + \textbf{V}_{j+1}\overline{\textbf{H}}_{j},
\]
leading to the cheap residual update given by 
\[
    \textbf{r}_{j} = \widehat{\textbf{r}} - (\textbf{A}\textbf{U}\textbf{B}_{j} + \textbf{V}_{j+1}\overline{\textbf{H}}_{j})\textbf{y}_{j} + \textbf{A}\textbf{U}\textbf{B}_{j}\textbf{y}_{j} \\ = \hat{\textbf{r}} - \textbf{V}_{j+1}\overline{\textbf{H}}_{j}.
\]

Finally, as with standard Krylov subspace methods, different types of augmented methods are distinguished based off the choice of constraint spaces $\widetilde{\mathcal{V}}_{j}$ and $\widetilde{\mathcal{U}}$. For example, by taking $\widetilde{\mathcal{U}} = \mathcal{C} := \textbf{A} \mathcal{U}$ and $
\widetilde{\mathcal{V}}_{j} = \textbf{A} \mathcal{K}_{j}((\textbf{I} - \Phi_{ \mathcal{C},\mathcal{C}^{\perp}})\textbf{A},(\textbf{I} - \Phi_{ \mathcal{C},\mathcal{C}^{\perp}})\textbf{r}_{0} ),
$
along with imposing $\textbf{C}$ to have orthonormal columns,
we recover the GCRO-DR algorithm \cite{parks2006recycling}.

\section{A framework for unprojected augmented Krylov Subspace Methods}\label{sec:alternative_formulation}
The work in this section constitutes the main contributions of this paper. We show that the projected problem (\ref{eq:framework_projectedproblem}) has an unprojected equivalent, and that viewing the projected problem in its unprojected form, allows the framework \cite{ gaul2013framework, gutknecht2014deflated, gaul2014recycling,gutknecht2012spectral, Soodhalter2020ASO} to easily adapt to describe the class of \textit{unprojected} augmented Krylov subspace methods. We use the framework to derive the first unprojected augmented and recycled Full Orthogonalization Method and demonstrate its effectiveness using numerical experiments. We then show how existing unprojected methods such as R$^{3}$GMRES and Augmented Conjugate Gradients fit within the framework. 

\vspace{0.4cm}
\begin{proposition}
\normalfont The projected problem (\ref{eq:framework_projectedproblem}) is  equivalent to the unprojected problem
\\ \newline
\noindent\fbox{%
    \parbox{\textwidth}{
  Find $\textbf{t}_{j} \in \mathcal{V}_{j}$  as an approximate solution to the problem 
\begin{equation}\label{eq:unprojected_problem}
    \textbf{A} \textbf{x} =  \textbf{r}_{0} 
\end{equation}
 subject to $\textbf{r}_{j} \perp (\textbf{I} - \Phi_{\textbf{A} \mathcal{U},\widetilde{\mathcal{U}}^{\perp}} )^{*}\widetilde{\mathcal{V}}_{j} := \widetilde{\mathcal{W}}_{j}$.
    }
}
\end{proposition}
\vspace{0.5cm}
\begin{proof} 
We denote the residual and solution of the projected problem (\ref{eq:framework_projectedproblem}) with the superscript $(1)$ and that of the unprojected problem (\ref{eq:unprojected_problem}) with a superscript (2). Application of the residual constraint of the projected problem yields
\begin{align*}
    \textbf{0} = \widetilde{\textbf{V}}^{*} \textbf{r}^{(1)}_{j} = \widetilde{\textbf{V}}^{*}_{j}((\textbf{I}-\Phi_{\textbf{A} \mathcal{U},\widetilde{\mathcal{U}}^{\perp}}) \textbf{r}_{0} - (\textbf{I}-\Phi_{\textbf{A} \mathcal{U},\widetilde{\mathcal{U}}^{\perp}}) \textbf{A} \textbf{V}_{j} \textbf{y}_{j}^{(1)} )
    \\ \implies \textbf{y}_{j}^{(1)} = (\widetilde{\textbf{V}}^{*}_{j} (\textbf{I}-\Phi_{\textbf{A} \mathcal{U},\widetilde{\mathcal{U}}^{\perp}}) \textbf{A} \textbf{V}_{j})^{-1} \widetilde{\textbf{V}}_{j}^{*}(\textbf{I}-\Phi_{\textbf{A} \mathcal{U},\widetilde{\mathcal{U}}^{\perp}}) \textbf{r}_{0},
\end{align*}
 and application of the residual constraint of the unprojected problem yields
\begin{align*}
 \textbf{0} = \widetilde{\textbf{V}}_{j}^{*}(\textbf{I} - \Phi_{\textbf{A} \mathcal{U},\widetilde{\mathcal{U}}^{\perp}} ) \textbf{r}_{j}^{(2)} = \widetilde{\textbf{V}}_{j}^{*}(\textbf{I} - \Phi_{\textbf{A} \mathcal{U},\widetilde{\mathcal{U}}^{\perp}} ) ( \textbf{r}_{0} - \textbf{A} \textbf{V}_{j} \textbf{y}_{j}^{(2)}) \\ \implies
 \textbf{y}_{j}^{(2)} = (\widetilde{\textbf{V}}^{*}_{j} (\textbf{I}-\Phi_{\textbf{A} \mathcal{U},\widetilde{\mathcal{U}}^{\perp}}) \textbf{A} \textbf{V}_{j})^{-1} \widetilde{\textbf{V}}_{j}^{*}(\textbf{I}-\Phi_{\textbf{A} \mathcal{U},\widetilde{\mathcal{U}}^{\perp}}) \textbf{r}_{0}.
\end{align*}
Thus $\textbf{y}^{(1)}_{j} = \textbf{y}_{j}^{(2)}$, and since both problems have the same search space $\mathcal{V}_{j}$ we have that in \textit{exact arithmetic}  $\textbf{x}_{j}^{(1)} = \textbf{x}_{j}^{(2)}$.
\end{proof}

For the unprojected form of the problem we compute  $\textbf{y}_{j}$ by solving the $j \times j$ linear system
\begin{equation} \label{projectedlinearsystem}
(\widetilde{\textbf{V}}_{j}^{*}( \textbf{I} - \Phi_{\textbf{A} \mathcal{U},\widetilde{\mathcal{U}}^{\perp}}) \textbf{A} \textbf{V}_{j}) \textbf{y}_{j} = \widetilde{\textbf{V}}^{*}_{j} ( \textbf{I} - \Phi_{\textbf{A} \mathcal{U},\widetilde{\mathcal{U}}^{\perp}}) \textbf{r}_{0}.
\end{equation}
The natural choice of subspace $\mathcal{V}_{j}$ for solving (\ref{projectedlinearsystem}) is the unprojected Krylov subspace $\mathcal{V}_{j} = \mathcal{K}_{j}(\textbf{A},\textbf{r}_{0})$.
The coefficient vector $\textbf{z}_{j}$ can be computed using (\ref{zcorrection}) in the same way as in the original formulation of the framework and is not dependent on whether the projected or unprojected problem is solved.

An important contrast between the projected and unprojected formulation can now be made. For any augmented method, in the projected formulation, problem (\ref{eq:framework_projectedproblem}) for $\textbf{y}_{j}$ is equivalent to applying the corresponding standard Krylov subspace method to an appropriately projected problem because the search and constraint space are $\mathcal{V}_{j}$ and $\widetilde{\mathcal{V}}_{j}$ respectively. An analogy of this result generally does not hold in the unprojected formulation, since in (\ref{eq:unprojected_problem}) the constraint space is no longer $\widetilde{\mathcal{V}}_{j}$, but is instead $\widetilde{\mathcal{W}}_{j} = (\textbf{I} - \Phi_{\textbf{A} \mathcal{U},\widetilde{\mathcal{U}}^{\perp}} )^{*}\widetilde{\mathcal{V}}_{j}$. In other words, the constraint space is composed of both a Krylov subspace and the augmentation subspace. Thus the problem for $\textbf{y}_{j}$ in the unprojected formulation (\ref{eq:unprojected_problem}) does not follow the same structure as the problem solved by the corresponding standard Krylov subspace method, unless additional constraints are imposed on the augmentation subspace. 

\vspace{0.4cm}
\begin{proposition}\label{prop1}\normalfont
For any unprojected augmented Krylov subspace method with search space $\mathcal{V}_{j} + \mathcal{U}$ and constraint space $\widetilde{\mathcal{V}}_{j} + \widetilde{\mathcal{U}}$,  in the special case 
when  $\textbf{A} \mathcal{U} \perp \widetilde{\mathcal{V}}_{j}$, then (\ref{eq:unprojected_problem}) is equivalent to applying the corresponding standard Krylov subspace method to an unprojected problem.
\end{proposition}
\begin{proof}
In order for (\ref{eq:unprojected_problem}) to be equivalent to applying the corresponding standard Krylov subspace method to an unprojected problem, we require $\widetilde{\mathcal{W}}_{j} =  \widetilde{\mathcal{V}}_{j}$. Since
\[
  \widetilde{\mathcal{W}}_{j} = (\textbf{I} - \Phi_{\textbf{A} \mathcal{U}, \widetilde{\mathcal{U}}^{\perp}})^{*} \widetilde{\mathcal{V}}_{j}
   = \widetilde{\mathcal{V}}_{j} - \widetilde{\mathcal{U}} (\widetilde{\mathcal{U}}^{*} \textbf{A} \mathcal{U})^{-*} (\textbf{A}\mathcal{U})^{*} \widetilde{\mathcal{V}}_{j},
 \]
this is satisfied when $
     \textbf{A} \mathcal{U} \perp \widetilde{\mathcal{V}}_{j}.$
\end{proof}

\vspace{0.4cm}
We now show how to derive new augmented unprojected Krylov subspace methods directly from the unprojected formulation of the framework by choosing appropriate choices of constraint spaces. We derive the first unprojected recycled Full Orthogonalization Method in the next subsection and in Subsection \ref{subsec:r3gmres} we show how the R$^{3}$GMRES algorithm can naturally be derived from the framework. In Subsection \ref{sec:augCG} we show how Proposition \ref{prop1} influences the derivation of the unprojected augmented Conjugate Gradient Algorithm \cite{erhel2000augmented}.

\subsection{Unprojected recycled FOM}\label{subsec:unprojFOM}
We develop the first unprojected recycled FOM using the appropriate choices of subspace for a FOM method, $\mathcal{V}_{j} = \widetilde{\mathcal{V}}_{j} = \mathcal{K}_{j}(\textbf{A},\textbf{b})$ and $\widetilde{\mathcal{U}} = \mathcal{U}$. We have the following proposition.
\begin{proposition} \normalfont
For the choice of subspaces $\mathcal{V}_{j} = \widetilde{\mathcal{V}}_{j} = \mathcal{K}_{j}(\textbf{A},\textbf{b})$, $\widetilde{\mathcal{U}} = \mathcal{U}$, if we define $\textbf{C} := \textbf{A} \textbf{U}$, then equation (\ref{projectedlinearsystem}) reduces to
\begin{equation}\label{FOMresult} \normalfont
    (\textbf{H}_{j} - \textbf{V}_{j}^{*} \textbf{C} (\textbf{U}^{*} \textbf{C})^{-1} \textbf{U}^{*}\textbf{V}_{j+1}\overline{\textbf{H}}_{j})\textbf{y}_{j} = \textbf{V}^{*}_{j}(\textbf{I} - \textbf{C}(\textbf{U}^{*} \textbf{C})^{-1} \textbf{U}^{*}) \textbf{r}_{0}.
\end{equation}
The coefficients vector $\textbf{z}_{j}$ is then
$\textbf{z}_{j} = (\textbf{U}^{*} \textbf{C})^{-1} \textbf{U}^{*} \textbf{r}_{0} - \textbf{B}_{j} \textbf{y}_{j}$.
\end{proposition}
\begin{proof} 
From the definitions of $\mathcal{U}$ and $\widetilde{\mathcal{U}}$ we have $\Phi_{\textbf{A} \mathcal{U}, \tilde{\mathcal{U}}^{\perp}} = \textbf{C}(\textbf{U}^{*} \textbf{C})^{-1} \textbf{U}^{*}$. 
Application of the Arnoldi process to build a basis for $\mathcal{K}_{j}(\textbf{A},\textbf{r}_{0})$ yields the usual Arnoldi relation (\ref{eq:arnoldi}), allowing the linear system (\ref{eq:probfory}) to be rewritten as
\[
    (\textbf{V}_{j}^{*}  \textbf{V}_{j+1} \overline{\textbf{H}}_{j} - \textbf{V}_{j}^{*} \textbf{C}(\textbf{U}^{*} \textbf{C})^{-1} \textbf{U}^{*} \textbf{V}_{j+1} \overline{\textbf{H}}_{j})\textbf{y}_{j} = \textbf{V}_{j}^{*}(\textbf{I} - \textbf{C}(\textbf{U}^{*} \textbf{C})^{-1} \textbf{U}^{*}) \textbf{r}_{0}.
\]
Orthonormality of $\textbf{V}_{j}$ then yields (\ref{FOMresult}), while the coefficients vector $\textbf{z}_{j}$ arises from (\ref{zcorrection}) with $\widetilde{\mathcal{U}} = \mathcal{U}$.
\end{proof}

Upon computing $\textbf{y}_{j}$ the full solution and residual approximations are given by
\begin{align*}
    \textbf{x}_{j} = \textbf{x}_{0} + \textbf{V}\textbf{y}_{j} - \textbf{U} \textbf{B}_{j} \textbf{y}_{j} + \textbf{U}(\textbf{U}^{*} \textbf{C})^{-1} \textbf{U}^{*}\textbf{r}_{0},
   \\ \textbf{r}_{j} = \textbf{r}_{0} - \textbf{V}_{j+1} \overline{\textbf{H}}_{j} \textbf{y}_{j} - \textbf{C}(\textbf{U}^{*} \textbf{C})^{-1} \textbf{U}^{*} \textbf{r}_{0} + \textbf{C} \textbf{B}_{j} \textbf{y}_{j}.
\end{align*}
The Ritz procedure for for updating the recycling subspace in unprojected recycled FOM can be derived by applying the following projection. 
\[
    \mbox{Find} \hspace{0.4cm} \textbf{y} \in  \mathcal{U} + \mathcal{K}_{j}(\textbf{A},\textbf{b}) \hspace{0.4cm} \mbox{such that} \hspace{0.4cm} \textbf{A} \textbf{y} - \mu \textbf{y} \perp  \mathcal{U} + \mathcal{K}_{j}(\textbf{A},\textbf{b}) \hspace{0.4cm} \mbox{for some} \hspace{0.1cm} \mu \in \mathbb{C},
\]
which leads to a small $(k + j) \times (k + j)$ eigenproblem for a vector $\textbf{g} \in \mathbb{C}^{k+j}$ given by
\begin{equation}\label{eq:aug_ritz_eig_prob_unprojrfom}
  \mu \begin{bmatrix}
      \textbf{U}^{*}\textbf{U} & \textbf{U}^{*}\textbf{V}_{j} \\
      \textbf{V}_{j}^{*} \textbf{U} & \textbf{I}
  \end{bmatrix} \textbf{g} = \begin{bmatrix}
      \textbf{U}^{*}\textbf{C} & \textbf{U}^{*}\textbf{V}_{j} \\ \textbf{V}_{j}^{*}\textbf{C} & \textbf{H}_{j}
  \end{bmatrix} \textbf{g}.
\end{equation}

Imposing orthonormality of the columns of $\textbf{U}$ allows the eigenproblem (\ref{eq:aug_ritz_eig_prob_unprojrfom}) to be further simplified. Upon solving (\ref{eq:aug_ritz_eig_prob_unprojrfom}) for $\textbf{g}$ we can compute $\textbf{y}$ via $\textbf{y} = \begin{bmatrix}
    \textbf{U} & \textbf{V}_{j}
\end{bmatrix}\textbf{g}$.

In Algorithm \ref{alg:unprojrFOM} below we outline a sketch of unprojected recycled FOM for solving a particular linear system in a sequence. For simplicity, our presentation assumes a previous linear system in the sequence has been solved and an orthonormal basis for a recycling subspace has been constructed and stored in $\textbf{U}$. It should be noted that this is not the case when solving the first system in the sequence, in which case a cycle of standard FOM should be applied and used to construct the first recycling subspace. After this cycle  unprojected recycled FOM can then be applied for subsequent cycles and system solves. 

\begin{algorithm}[H]
\caption{Unprojected Recycled FOM \label{alg:unprojrFOM}}
\begin{algorithmic}[1]
\STATE{\textbf{Input:} $\textbf{A} \in \mathbb{C}^{n \times n}$ , $\textbf{x}, \textbf{b} \in \mathbb{C}^{n}$, $\textbf{U} \in \mathbb{C}^{n \times k}$ basis for $\mathcal{U}$, and $\textbf{C}$ such that $\textbf{C} = \textbf{A} \textbf{U}$}
\STATE{$\textbf{r}_{0} \leftarrow \textbf{b} - \textbf{A} \textbf{x}$}
\WHILE{ not converged }
\STATE{ Build $\mathcal{K}_{j}(\textbf{A},\textbf{r}_{0})$ via the Arnoldi algorithm to generate $\textbf{V}_{j+1}$ and $\overline{\textbf{H}}_{j}$}
\STATE{Define $\textbf{B}_{j} := (\textbf{U}^{*}  \textbf{C})^{-1} \textbf{U}^{*} \textbf{V}_{j+1} \overline{\textbf{H}}_{j}$}
\STATE{$\textbf{r} \leftarrow \textbf{r}_{0} - \textbf{C}(\textbf{U}^{*} \textbf{C})^{-1} \textbf{U}^{*} \textbf{r}_{0}$}
\STATE{ Solve  $(\textbf{H}_{j} - \textbf{V}_{j}^{*} \textbf{C} (\textbf{U}^{*} \textbf{C})^{-1} \textbf{U}^{*}\textbf{V}_{j+1}\overline{\textbf{H}}_{j})\textbf{y}_{j} = \textbf{V}^{*}_{j}\textbf{r}$}
\STATE{ $\textbf{x} \leftarrow \textbf{x} + \textbf{V}_{j}\textbf{y}_{j} - \textbf{U}\textbf{B}_{j}\textbf{y}_{j} + \textbf{U}(\textbf{U}^{*} \textbf{C})^{-1} \textbf{U}^{*} \textbf{r}_{0}$}
\STATE{ $\textbf{r}_{0} \leftarrow \textbf{r} - \textbf{V}_{j+1} \overline{\textbf{H}}_{j} \textbf{y}_{j} + \textbf{C} \textbf{B}_{j} \textbf{y}_{j}$}
\STATE{Update $\mathcal{U}$ and $\mathcal{C}$ for next cycle}
\ENDWHILE
\STATE{Update $\mathcal{U}$ and $\mathcal{C}$ for next system}
\end{algorithmic}
\end{algorithm}

We have implemented unprojected recycled FOM into MATLAB to demonstrate its behaviour as a recycling method. We compare our results to standard FOM and the recycled projected FOM derived from the projected formulation of the framework.

The experiment was run on the built in MATLAB Neumann matrix of size $22,500 \times 22,500$. A sequence of $5$ linear systems was solved with fixed Neumann matrix, and changing random right hand sides. An Arnoldi cycle length of $j = 90$ and a recycling subspace dimension of $k = 20$ was used. For the unprojected method the Ritz recycling procedure outlined in this section was used to construct the recycling subspace. A similar Ritz procedure was also used to construct the recycling subspace for projected recycled FOM.

\begin{figure}[H]
   \centering
    \includegraphics[width=8cm]{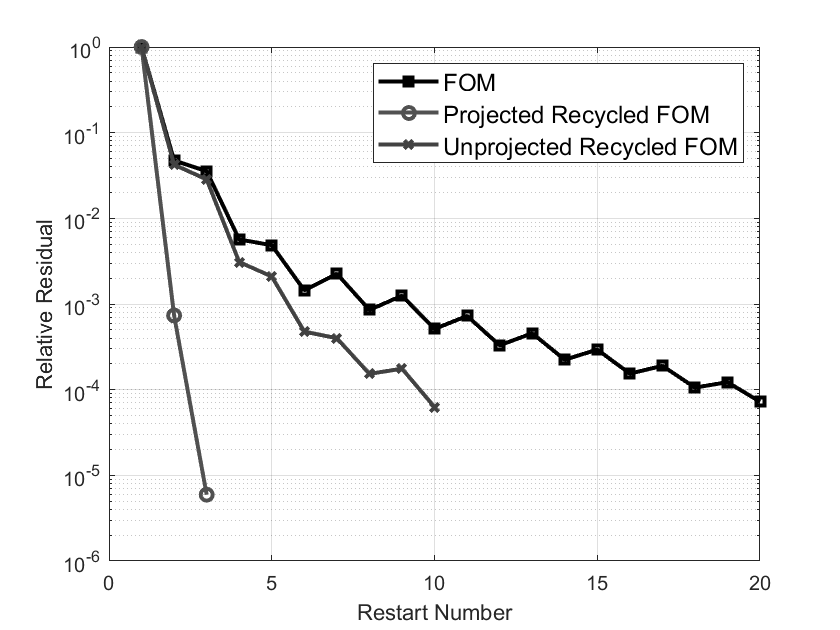}
   \caption{Convergence curves for the final linear system in a sequence of $5$ linear systems with fixed Neumann matrix of size $22,500 \times 22,500$ and changing random right hand sides.} 
  \label{fig:convergence_unprojectedrecycledfom}
\end{figure}

It is well known that projected recycling methods typically converge faster than unprojected methods and this is clearly seen in the convergence plot in Figure \ref{fig:convergence_unprojectedrecycledfom}. We record the number of times the matrix $\textbf{A}$ needs to be applied to reach convergence for each system in the table below.

\begin{tabular}{ |p{2.3cm}||p{1cm}|p{2.5cm}|p{3cm}|  }
 \hline
 \multicolumn{4}{|c|}{\# MAT-VEC's} \\
 \hline
 System Number & FOM & Projected rFOM & Unprojected rFOM\\
 \hline
 1  &  1890   & 660  & 1650  \\
 2  &  3150   & 440 & 2530   \\
 3  &  2790   & 220 & 2200   \\
 4  &  2700   & 220 & 1870   \\
 5  &  1710   & 220 &  990  \\
 \hline
\end{tabular}

\subsection{R$^{3}$GMRES}
\label{subsec:r3gmres}
In this section we show how the R$^{3}$GMRES algorithm \cite{dong2014r3gmres} can easily be derived from the unprojected formulation of the framework presented in this paper. Similar to the GCRO-DR algorithm \cite{parks2006recycling}, the  R$^{3}$GMRES algorithm can be formulated in terms of a minimization problem  over an augmented Krylov subspace. In \cite[Theorem 4.1]{soodhalter2021note}, it was shown that the normal equations associated to this minimization problem
can be decoupled into a smaller problem over a Krylov subspace, followed by computing an additional correction from the augmentation subspace. We show how the alternative formulation of R$^{3}$GMRES derived in
\cite[Theorem 4.1]{soodhalter2021note} can also be derived more easily and naturally as a direct result of the unprojected formulation of the framework presented in this paper. We do this by showing that through appropriate choices of subspaces $\widetilde{\mathcal{V}}_{j}, \widetilde{\mathcal{U}}$ and $\mathcal{U}$ corresponding to a GMRES method, namely $\mathcal{V}_{j} = \mathcal{K}_{j}(\textbf{A},\textbf{r}_{0})$, $\widetilde{\mathcal{V}}_{j} = \textbf{A} \mathcal{V}_{j}$ and $\widetilde{\mathcal{U}} = \textbf{A}\mathcal{U}$, the linear system (\ref{projectedlinearsystem}) reduces to the same problem derived in \cite[Theorem 4.1]{soodhalter2021note}.
\begin{theorem}\label{thm:r3gmresthm} \normalfont
For the choices of subspaces  $\mathcal{V}_{j} = \mathcal{K}_{j}(\textbf{A},\textbf{r}_{0})$, $\widetilde{\mathcal{V}}_{j} = \textbf{A} \mathcal{V}_{j}$ and $\widetilde{\mathcal{U}} = \textbf{A}\mathcal{U}$ with orthonormal columns of the matrix
$\textbf{C} := \textbf{A}\mathcal{U}$, the problem (\ref{eq:probfory}) reduces to the following $j \times j$ linear system for $\textbf{y}_{j}$
\begin{equation}\label{r3gmreslinearsystem}
    (\overline{\textbf{H}}_{j}^{*} \overline{\textbf{H}}_{j} - \overline{\textbf{H}}_{j}^{*} \textbf{V}^{*}_{j+1} \textbf{C} \textbf{C}^{*} \textbf{V}_{j+1} \overline{\textbf{H}}_{j} ) \textbf{y}_{j} = \overline{\textbf{H}}_{j}^{*} \textbf{V}_{j+1}^{*}(\textbf{I} - \textbf{C} \textbf{C}^{*}) \textbf{r}_{0}.
\end{equation}
The coefficients vector  $\textbf{z}_{j}$ is then $ \textbf{z}_{j} = \textbf{C}^{*} \textbf{r}_{0} - \textbf{B}_{j} \textbf{y}_{j}$.
\end{theorem}

\begin{proof}
From the definitions of $\widetilde{\mathcal{U}}$ and $\mathcal{U}$ and, due to orthonormality of $\textbf{C}$ we have $\Phi_{\textbf{A} \mathcal{U}, \widetilde{\mathcal{U}}^{\perp}} = \Phi_{\mathcal{C},\mathcal{C}^{\perp}} = \textbf{C} \textbf{C}^{*}$. Application of the Arnoldi process to build a basis for $\mathcal{K}_{j}(\textbf{A},\textbf{r}_{0})$ yields the usual Arnoldi relation (\ref{eq:arnoldi}), allowing the linear system (\ref{eq:probfory}) to be rewritten as
\begin{align*}
 (\widetilde{\textbf{V}}_{j}^{*}( \textbf{I} - \textbf{C} \textbf{C}^{*}) \textbf{A} \textbf{V}_{j}) \textbf{y}_{j} = \widetilde{\textbf{V}}^{*}_{j} ( \textbf{I} - \textbf{C} \textbf{C}^{*}) \textbf{r}_{0} \\\implies
((\textbf{A}\textbf{V}_{j})^{*} \textbf{A} \textbf{V}_{j}  -  (\textbf{A}\textbf{V}_{j})^{*} \textbf{C} \textbf{C}^{*} \textbf{A} \textbf{V}_{j}) \textbf{y}_{j} = ( \textbf{A}\textbf{V}_{j})^{*} ( \textbf{I} - \textbf{C} \textbf{C}^{*}) \textbf{r}_{0}
\\ \implies ( (\textbf{V}_{j+1} \overline{\textbf{H}}_{j})^{*}(\textbf{V}_{j+1} \overline{\textbf{H}}_{j}) - (\textbf{V}_{j+1} \overline{\textbf{H}}_{j})^{*} \textbf{C} \textbf{C}^{*} \textbf{V}_{j+1} \overline{\textbf{H}}_{j}  ) \textbf{y}_{j} = (\textbf{V}_{j+1} \overline{\textbf{H}}_{j})^{*}( \textbf{I} - \textbf{C} \textbf{C}^{*}) \textbf{r}_{0}
\\ \implies (\overline{\textbf{H}}^{*}_{j} \overline{\textbf{H}}_{j} - \overline{\textbf{H}}_{j}^{*} \textbf{V}^{*}_{j+1} \textbf{C} \textbf{C}^{*} \textbf{V}_{j+1} \overline{\textbf{H}}_{j}) \textbf{y}_{j} = \overline{\textbf{H}}^{*}_{j} \textbf{V}_{j+1}^{*} ( \textbf{I} - \textbf{C} \textbf{C}^{*}) \textbf{r}_{0}.
\end{align*}
The coefficients vector $\textbf{z}_{j}$ arises from (\ref{zcorrection}) with $\widetilde{\mathcal{U}} = \textbf{A} \mathcal{U}$.
\end{proof}

The result of Theorem \ref{thm:r3gmresthm} is equivalent to that shown in \cite[Theorem 4.1]{soodhalter2021note}, but derived from the unprojected formulation of the framework outlined in this paper. We also note that upon solving for $\textbf{y}_{j}$, the 
solution and the residual updates for R$^{3}$GMRES are
\begin{align*}
    \textbf{x}_{j} = \textbf{x}_{0} + \textbf{V}_{j} \textbf{y}_{j} + \textbf{U} \textbf{C}^{*} \textbf{r}_{0} - \textbf{U} \textbf{B}_{j} \textbf{y}_{j}\\
   \textbf{r}_{j} = \textbf{r}_{0} - \textbf{V}_{j+1} \overline{\textbf{H}}_{j} \textbf{y}_{j} - \textbf{C} \textbf{C}^{*} \textbf{r}_{0} + \textbf{C} \textbf{B}_{j} \textbf{y}_{j}.
\end{align*}
We refer the reader to \cite{soodhalter2021note} for an overview of the implementation details of this formulation of R$^{3}$GMRES. An efficient MATLAB implementation is also provided \cite{kirk_m_soodhalter_2021_4975990}.

\subsection{Unprojected Augmented Short Recurrence Methods}\label{sec:augCG}
In the projected formulation of the framework, it has been shown in \cite{gutknecht2011deflated} (see also \cite{Soodhalter2020ASO}) that for a Symmetric Positive Definite (SPD) matrix $\textbf{A}$ the projected matrix $(\textbf{I} - \Phi_{\textbf{A} \mathcal{U}, \widetilde{\mathcal{U}}^{\perp}} ) \textbf{A}$ is also SPD and a standard Conjugate Gradient algorithm can be applied to the projected problem. It was shown in \cite{gutknecht2011deflated} that for an augmented MINRES for Hermitian indefinite systems, although the resulting projected matrix is non-Hermitian the Krylov subspace built from this matrix and the projected residual is equivalent to that built from a Hermitian matrix and so the usual short recurrence MINRES algorithm can also be applied in this setting. 

In the unprojected formulation however, as a a consequence of Proposition \ref{prop1} the SPD properties of the operator are generally not inherited into the linear system (\ref{projectedlinearsystem}) unless the condition $\textbf{A} \mathcal{U} \perp \widetilde{\mathcal{V}}_{j}$ is satisfied, leading to difficulties when developing short recurrence methods for SPD problems.

We discuss the Augmented Conjugate Gradient (AugCG) algorithm proposed in \cite{erhel2000augmented} and show how the method fits within the unprojected formulation of the framework while satisfying the conditions of proposition \ref{prop1} .

AugCG can be used to solve consecutive SPD linear systems of the form
\[
    \textbf{A}\textbf{x}_{1} = \textbf{b}_{1} \hspace{1cm} \textbf{A}\textbf{x}_{2} = \textbf{b}_{2}, \hspace{1cm} \textbf{b}_{1}, \textbf{b}_{2} \in \mathbb{C}^{n}.
\]

If the standard Conjugate Gradient algorithm is applied to the first system to generate the set of residuals residuals $\textbf{R}^{(1)}_{j} = \begin{bmatrix} \textbf{r}_{0}^{(1)}, \ldots , \textbf{r}_{j}^{(1)} \end{bmatrix}$ and search directions $\textbf{P}^{(1)}_{j} = \begin{bmatrix}\textbf{p}_{0}^{(1)}, \ldots ,\textbf{p}_{j}^{(1)}\end{bmatrix}$ such that 
\[
\mbox{span}(\textbf{R}^{(1)}_{j}) = \mbox{span}(\textbf{P}^{(1)}_{j}) = \mathcal{K}_{j}(\textbf{A},\textbf{r}_{0}^{(1)}),
\]
then the augmented Conjugate Gradient algorithm (AugCG) can  be applied to the second system by taking $\mathcal{U} = \mathcal{K}_{j}(\textbf{A},\textbf{r}_{0}^{(1)})$ and $\mathcal{V}_{j} = \mathcal{K}_{j}(\textbf{A}, \textbf{r}_{0}^{(2)})$ where $\textbf{r}_{0}^{(2)}$ is the starting initial residual used to solve system $2$, chosen such that $\textbf{r}_{0}^{(2)} \perp \mathcal{U}$. At iteration $j+1$ we impose that
 the residuals $\textbf{r}_{j+1}^{(2)}$ be orthogonal to $\mathcal{V}_{j}$ and the descent descent directions $\textbf{p}^{(2)}_{j+1}$ are $\textbf{A}$-conjugate to both $\mathcal{U}$ and $\mathcal{V}_{j}$. The starting  descent direction $\textbf{p}_{0}^{(2)}$ is thus computed via
\begin{equation}\label{eq:AugCG_conditions}
    \textbf{p}_{0}^{(2)} = (\textbf{I} - \textbf{U}(\textbf{U}^{*} \textbf{A} \textbf{U})^{-1} (\textbf{A} \textbf{U})^{*}) \textbf{r}_{0}^{(2)}.
\end{equation}
This $\textbf{A}$-conjugacy condition is maintained by line 11 in Algorithm \ref{alg:augCG} below.
By the end of the process, it was shown in \cite{erhel2000augmented} using induction that the relations 
\begin{equation}
\textbf{r}_{j+1}^{(2)} \perp \mathcal{K}_{j}(\textbf{A}, \textbf{r}_{0}^{(2)}) \hspace{1cm} \textbf{U}^{*}\textbf{A} \textbf{p}_{j+1}^{(2)} = \textbf{0},
\end{equation}
are satisfied.

The AugCG algorithm as outlined in \cite{erhel2000augmented} is outlined below.
\begin{algorithm}[H]
\caption{Augmented Conjugate Gradient (AugCG) \label{alg:augCG}}
\begin{algorithmic}[1]
\STATE{ \textbf{Input:} $\textbf{A} \in \mathbb{C}^{n \times n}$ , $ \textbf{b} \in \mathbb{C}^{n}$, initial residual $\textbf{r}_{0} \in \mathbb{C}^{n}$ from initial approximation $\textbf{x}_{0} \in \mathbb{C}^{n}$,  basis of Krylov subspace for previous system $\textbf{U}$, matrix $\textbf{C} = \textbf{A}\textbf{U}$.}
\STATE{ $\textbf{r}_{0} = \textbf{r}_{0} - \textbf{C}(\textbf{U}^{*}\textbf{C})^{-1} \textbf{U}^{*}\textbf{r}_{0}$}
\STATE{$\textbf{x}_{0} = \textbf{x}_{0} + \textbf{U}(\textbf{U}^{*} \textbf{C})^{-1} \textbf{U}^{*} \textbf{r}_{0}$}
\STATE{$\textbf{p}_{0} = \textbf{r}_{0} - \textbf{U}(\textbf{U}^{*}\textbf{C})^{-1} \textbf{C}^{*} \textbf{r}_{0}$}
\FOR{$j=0,1, \hdots$ until convergence}
 \STATE{$\alpha_{j} = \textbf{r}_{j}^{T} \textbf{r}_{j} / \textbf{p}_{j}^{T} \textbf{A} \textbf{p}_{j}$}
 \STATE{$\textbf{x}_{j+1} = \textbf{x}_{j} + \alpha_{j} \textbf{p}_{j}$}
 \STATE{$\textbf{r}_{j+1} = \textbf{r}_{j} - \alpha_{j} \textbf{A} \textbf{p}_{j}$}
 \STATE{ $\beta_{j} = \textbf{r}_{j+1}^{*} \textbf{r}_{j+1}/\textbf{r}_{j}^{*}\textbf{r}_{j}$}
\STATE{$\mu_{k+1} = \textbf{w}_{k}^{*} \textbf{A} \textbf{r}_{j+1} / \textbf{w}_{k}^{T} \textbf{A} \textbf{w}_{k}$}
 \STATE{$\textbf{p}_{j+1} = \textbf{r}_{j+1} + \beta_{j} \textbf{p}_{j} - \mu_{k+1} \textbf{w}_{k}$}
\ENDFOR
\end{algorithmic}
\end{algorithm}

We now derive the AugCG algorithm directly from the framework presented in this paper and demonstrate the implications Proposition  \ref{prop1} has on the development of AugCG and discuss why it is necessary to take $\mathcal{U} = \mathcal{K}_{j}(\textbf{A}, \textbf{r}_{0}^{(1)})$.

We assume $\textbf{A}$ is symmetric positive definite and that $\widetilde{\mathcal{U}} = \mathcal{U}$ and $\widetilde{\mathcal{V}}_{j} = \mathcal{V}_{j} = \mathcal{K}_{j}(\textbf{A},\textbf{r}_{0})$ in which case (\ref{eq:eqfory}) reduces to 
\begin{equation}\label{eq:CGprobfory}
    (\textbf{T}_{j} - \textbf{V}_{j}^{*} \textbf{C} (\textbf{U}^{*} \textbf{C})^{-1} \textbf{U}^{*}\textbf{A} \textbf{V}_{j})\textbf{y}_{j} = \textbf{V}^{*}_{j}(\textbf{I} - \textbf{C}(\textbf{U}^{*} \textbf{C})^{-1} \textbf{U}^{*}) \textbf{r}_{0}.
\end{equation}
The solution and residual updates are
\begin{align*}
    \textbf{x}_{j} = \textbf{x}_{0} + \textbf{V}\textbf{y}_{j} - \textbf{U} \textbf{B}_{j} \textbf{y}_{j} + \textbf{U}(\textbf{U}^{*} \textbf{C})^{-1} \textbf{U}^{*}\textbf{r}_{0},
    \\ \textbf{r}_{j} = \textbf{r}_{0} - \textbf{V}_{j+1} \overline{\textbf{T}}_{j} \textbf{y}_{j} - \textbf{C}(\textbf{U}^{*} \textbf{C})^{-1} \textbf{U}^{*} \textbf{r}_{0} + \textbf{C} \textbf{B}_{j} \textbf{y}_{j},
\end{align*}
with $\textbf{B}_{j} = (\textbf{U}^{*}\textbf{A} \textbf{U})^{-1} \textbf{U}^{*}\textbf{A}\textbf{V}_{j}$.

A practical implementation of a short recurrence method requires the SPD structure of the matrix $\textbf{A}$ to be inherited into the problem for $\textbf{y}_{j}$. While this was shown to be true for the projected formulation of augmented CG in \cite{Soodhalter2020ASO}, it is not generally true for in the unprojected case since the coefficient matrix of the linear system (\ref{eq:CGprobfory}) is generally not SPD. This is a direct consequence of the fact that the condition $\textbf{A} \mathcal{U} \perp \widetilde{\mathcal{V}}_{j}$ from proposition \ref{prop1} generally does not hold. In order to ensure this condition holds we impose (\ref{eq:AugCG_conditions}) which leads to 
\[ \textbf{U}^{*}\textbf{A}\textbf{V}_{j} = \textbf{0} \implies \textbf{V}_{j}^{*} \textbf{C} (\textbf{U}^{*} \textbf{C})^{-1} \textbf{U}^{*}\textbf{A} \textbf{V}_{j} = \textbf{0},  \hspace{0.5cm}   \textbf{B}_{j} = \textbf{0}
\]
This ensures the coefficient matrix of the problem (\ref{eq:CGprobfory}) is now SPD and a standard CG iteration can now be applied to (\ref{eq:CGprobfory}) in line with Proposition \ref{prop1}.
Imposing the conditions (\ref{eq:AugCG_conditions}) can only be done using the fact that the residuals of the standard FOM iteration and multiples of the Arnoldi vectors, imposing a strict choice of augmentation subspace $\mathcal{U}$ to be the Krylov subspace from the previous system solve.

\section{Conclusion}\label{sec:conclusions}
In this paper we have taken an alternative view of a general framework describing projected augmented Krylov subspace methods. We have shown that by viewing the projected problem arising in this framework in terms of its unprojected equivalent, the framework  easily and naturally extends to describe unprojected augmented Krylov subspace methods. We demonstrated this by using the framework to derive a new unprojected recycled FOM and its associated Ritz problem, as well as existing methods such as R$^{3}$GMRES.

It is well known from the original formulation of the framework that a projected Krylov subspace method involves applying the corresponding standard Krylov subspace method to an appropriately projected problem. We showed in this work that an analogy of this result does not hold for unprojected augmented methods unless additional assumptions are imposed on the augmentation subspace. This leads to practical limitations when developing unprojected augmented short recurrence methods for Hermitian problems, but it is still possible to develop such methods from the framework under certain assumptions on the recycling subspace. We demonstrated this using Augmented Conjugate Gradient (AugCG) as an example.

\section{Acknowledgements}
This work was jointly funded by the Irish Research Council Government of Ireland Postgraduate Scholarship Programme, The Hamilton Scholars, and a Trinity College Dublin postgraduate award.

\printbibliography %Prints bibliography
\end{document}